\documentclass[10pt]{article}
\font\smallit=cmti10
\font\smalltt=cmtt10

\usepackage{amssymb,latexsym,amsmath,epsfig,amsthm,hyperref} %% Add other packages as necessary

\makeatletter

\renewcommand\section{\@startsection {section}{1}{\z@}
{-30pt \@plus -1ex \@minus -.2ex}
{2.3ex \@plus.2ex}
{\normalfont\normalsize\bfseries\boldmath}}

\renewcommand\subsection{\@startsection{subsection}{2}{\z@}
{-3.25ex\@plus -1ex \@minus -.2ex}
{1.5ex \@plus .2ex}
{\normalfont\normalsize\bfseries\boldmath}}

\renewcommand{\@seccntformat}[1]{\csname the#1\endcsname. }

\makeatother

\newcommand{\bburl}[1]{\textcolor{blue}{\url{#1}}}
\newcommand{\seqnum}[1]{\href{https://oeis.org/#1}{\rm \underline{#1}}}

\newtheorem{theorem}{Theorem}
\newtheorem{lemma}{Lemma}
\newtheorem{proposition}{Proposition}
\newtheorem{corollary}{Corollary}

\theoremstyle{definition}

\newtheorem{remark}{Remark}

\begin{document}

\begin{center}
\uppercase{\bf Schreier Multisets and the $\mathbf S$-step Fibonacci Sequences}
\vskip 20pt
{\bf H\`ung Vi\d{\^e}t Chu}\\
{\smallit Department of Mathematics, University of Illinois Urbana-Champaign, Urbana, IL, USA}\\
{\tt hungchu2@illinois.edu}\\ 
\vskip 10pt
{\bf Nurettin Irmak}\\
{\smallit Department of Engineering Basic Sciences, Konya Technical University, Konya, T\"urkiye}\\
{\tt irmaknurettin@gmail.com; \tt nirmak@ktun.edu.tr}\\ 
\vskip 10pt
{\bf Steven J. Miller}\\
{\smallit Department of Mathematics and Statistics, Williams College, Williamstown, MA, USA}\\
{\tt sjm1@williams.edu; \tt Steven.Miller.MC.96@aya.yale.edu}\\ 
\vskip 10pt
{\bf L\'{a}szl\'{o} Szalay}\\
{\smallit Department of Mathematics, J. Selye University, Kom\'{a}rno, Slovakia}\\
{\smallit Institute of Informatics and Mathematics, University of Sopron, Sopron, Hungary}
{\tt szalay.laszlo@uni-sopron.hu}\\ 
\vskip 10pt
{\bf Sindy Xin Zhang}\\
{\smallit Department of Mathematics, University at Buffalo, Buffalo, NY, USA}\\
{\smallit Department of Mathematics, Tufts University, Medford, MA, USA}\\
{\tt xzhang99@buffalo.edu; \tt sindy.zhang@tufts.edu}\\ 
\end{center}
\vskip 20pt
\centerline{\smallit Received: , Revised: , Accepted: , Published: } % We will fill in the dates
\vskip 30pt

\centerline{\bf Abstract}
\noindent
Inspired by the surprising relationship (due to A. Bird) between Schreier sets and the Fibonacci sequence,  we introduce Schreier multisets and connect these multisets with the $s$-step Fibonacci sequences, defined, for each $s\geqslant 2$, as: $F^{(s)}_{2-s} = \cdots = F^{(s)}_0 = 0$, $F^{(s)}_1 = 1$, and 
$F^{(s)}_{n}  = F^{(s)}_{n-1} + \cdots + F^{(s)}_{n-s}, \mbox{ for } n\geqslant 2$.
Next, we use Schreier-type conditions on multisets to retrieve a family of sequences which satisfy a recurrence of the form $a(n) = a(n-1) + a(n-u)$, with $a(n) = 1$ for $n = 1,\ldots, u$. Finally, we study nonlinear Schreier conditions and show that these conditions are related to integer decompositions, each part of which is greater than the number of parts raised to some power. 

\pagestyle{myheadings}
\markright{\smalltt INTEGERS: 23 (2023)\hfill}
\thispagestyle{empty}
\baselineskip=12.875pt
\vskip 30pt

\section{Introduction}
A set $A\subset\mathbb{N}$ is said to be \textit{Schreier} if $\min A \geqslant |A|$. These sets have been extensively studied both in Banach space theory and in Ramsey theory. In 2012, Bird  \cite{B} showed that 
$$
|\{A\subset\{1,2,\dots,n\}\, :\, n\in A\mbox{ and }\min A\geqslant |A|\}|\ =\ F_n, \mbox{ for all } n\geqslant 1,
$$ 
where $F_n$ is the $n$\textsuperscript{th} Fibonacci number defined as $F_0 = 0$, $F_1 = 1$, and $F_n = F_{n-1}+F_{n-2}$ for $n\geqslant 2$. Since then, various other recurrence relations have been established by Schreier-type conditions: \cite{BCF, BGHH, C1, CMX}. In this paper, we shall use the Schreier condition to obtain the $s$-step Fibonacci (or the $s$-Fibonacci) sequences, thus generalizing Bird's result. Particularly, fix $s\geqslant 2$ and define the $s$-step Fibonacci sequence as: $F^{(s)}_{2-s} = \cdots = F^{(s)}_0 = 0$, $F^{(s)}_1 = 1$, and 
$$F^{(s)}_{n}  \ =\ F^{(s)}_{n-1} + \cdots + F^{(s)}_{n-s}, \mbox{ for } n\geqslant 2.$$
For each $n\in\mathbb{N}$, set
$$\mathcal{A}^{(s-1)}_{n} \ =\ \{A\subset \{\underbrace{1, \ldots, 1}_{s-1}, \ldots, \underbrace{n-1, \ldots, n-1}_{s-1}, n\}\,:\, n\in A\mbox{ and }\min A\geqslant |A|\}.$$

We state our first new result.
\begin{theorem}\label{m1}
For $n\in\mathbb{N}$ and $s\geqslant 2$, it holds that $|\mathcal{A}^{(s-1)}_{n}| = F^{(s)}_n$.
\end{theorem}
Note that in the case $s = 2$, Theorem \ref{m1} gives Bird's result. 

The first named author of the present paper \cite{C3} recently discovered another way to generate the Fibonacci sequence from Schreier-type sets. Specifically, for $n\geqslant 1$, 
$$ |\mathcal{D}_n \ :=\ \{A\subset \{1, \ldots, n\}\,:\, \mbox{ either }A = \emptyset\mbox{ or }(\max A - 1\in A\mbox{ and }\min A\geqslant |A|)\}|$$
is equal to $F_n$. 
This way does not fix the maximum of sets as what Bird did; however, once the maximum $a$ of a set is chosen, then $a-1$ must also be in the set. (This requirement is used to prove an injective map between $\mathcal{D}_{n+1}\backslash \mathcal{D}_{n}$ and $\mathcal{D}_{n-1}$.) Using multisets, we show yet another way to generate the Fibonacci sequence. Fix a sequence $\mathbf s = (s_n)_{n=1}^\infty$ of nonnegative integers satisfying $s_{n}\geqslant k$ for all $n\geqslant 2k+1$ and $k\geqslant 1$. Define
$$\mathcal{B}_n^{\mathbf s}\ :=\ \{B\subset \{\underbrace{1, \ldots, 1}_{s_1}, \underbrace{2, \ldots, 2}_{s_2}, \ldots, \underbrace{n, \ldots, n}_{s_n}\}\,:\, B = \emptyset \mbox{ or } \min B\geqslant 2|B|+1\}.$$
It turns out that $|\mathcal{B}^{\mathbf s}_{n}| = F_n$ for $n\in \mathbb{N}$. Indeed, we prove a more general result. Fix $u\geqslant 2$ and define the sequence $(K^{(u)}_n)_{n=1}^\infty$ as follows:
$$K^{(u)}_1 \ =\ \cdots\ \ =\ K^{(u)}_u \ = \ 1 \mbox{ and }K^{(u)}_{n}\ =\ K^{(u)}_{n-1} + K^{(u)}_{n-u}, \quad n\geqslant u+1.$$
Given a sequence $\mathbf s = (s_n)_{n=1}^\infty$ of nonnegative numbers, let
$$\mathcal{B}^{\mathbf s, u}_n\ :=\ \{B\subset \{\underbrace{1, \ldots, 1}_{s_1}, \underbrace{2, \ldots, 2}_{s_2}, \ldots, \underbrace{n, \ldots, n}_{s_n}\}\,:\, \min B\geqslant u|B|+1\}.$$

\begin{theorem}\label{m2}
Fix $u\geqslant 2$ and a sequence $\mathbf s = (s_n)_{n=1}^{\infty}$ of nonnegative integers such that 
$s_{n}\geqslant k$ for all $n\geqslant uk+1$ and $k\geqslant 1$.
Then 
$$|\mathcal{B}^{\mathbf s, u}_{n}| \ =\  K^{(u)}_n,\quad n\in \mathbb{N}.$$
\end{theorem}

Our next result is also related to the sequence $(K^{(u)}_n)_{n=1}^\infty$. We use colored multisets instead of uncolored multisets as in defining $\mathcal{A}^{(s-1)}_{n}$. In particular, for a fixed $s\geqslant 2$, assume that the positive integers $1,2\dots,n-1$ have colors from a palette of $(s-1)\geqslant 1$ tints. We denote by $i_j$ if the integer $i$ possesses color $j$. Let
$$
H^{(s-1)}_n\ :=\ \{1_1,1_2,\dots,1_{s-1},2_1,2_2,\dots,2_{s-1},\dots,(n-1)_1,(n-1)_2,\dots,(n-1)_{s-1},n\},
$$
where the number $n$ is uncolored. Put	
$$
\mathcal{C}^{(s-1)}_n \ :=\ \{C\subset H^{(s-1)}_n\, :\, n\in C\mbox{ and }\min C\geqslant |C|\}.
$$ 
\begin{theorem}\label{m4}
	For $n\in\mathbb{N}$ and $s\geqslant 2$, we have $$|\mathcal{C}^{(s-1)}_n|\ =\ K^{(s)}_{(s-1)(n-1)+1}.$$
\end{theorem}
We remark that one can give a direct recurrence relation of order $s$ between the terms of the sequence $\tau_{n}=K^{(s)}_{(s-1)(n-1)+1}$. Consider $s=4$, for example.
The sequence $(K^{(s)}_n)_{n=1}^\infty$ has values 
$$
\frame{1},1,1,\frame{1},2,3,\frame{4},5,7,\frame{10},14,19,\frame{26},36,50,\frame{69},95,131,\frame{181},250,345,\frame{476},\dots,
$$
and the framed entries satisfy the (direct) recursive relation 
$$\tau_n\ =\ \tau_{n-1}+3\tau_{n-2}+3\tau_{n-3}+\tau_{n-4},$$
see \seqnum{A003269} and \seqnum{A099234} in \cite{Sl}.

Our final results involve counting sets $A$ under nonlinear Schreier conditions of the form $\sqrt[s]{\min A} \geqslant |A|$ for some fixed $s\geqslant 2$. Such sets recently appeared in \cite[Proposition 6.10]{BDKOW}. For $n\in\mathbb{N}$ and  $p\in \mathbb{Z}_{\geqslant 0}$, let $K_{n, p}$ count the number of decompositions of $n$, where the smallest part is strictly greater than the number of parts raised to the $p$\textsuperscript{th} power. In notation,
$$K_{n,p}\ :=\ \left|\left\{(x_1, \ldots, x_k)\,:\, \sum_{i=1}^k x_i = n\mbox{ and } \min x_i > k^{p}\right\}\right|.$$
When $p = 0$, the sequence $K_{n,0} = F_{n-1}$ for $n\geqslant 1$ (see Corollary \ref{c1}.) When $p = 1$, the sequence $(K_{n,1})_{n=1}^\infty$ is ~\seqnum{A098132}. The first few values are $$0, 1, 1, 1, 1, 2, 3, 4, 5, 6, 7, 9, 12, 16, 21, 27, 34, 42, \ldots.$$
For example, $K_{13,1} = 12$ because we can write $13$ as
\begin{align*}
13&\ =\  3+10 \ =\ 10 + 3 \ =\  4 + 9 \ =\ 9 + 4 \ =\ 5 + 8 \ =\ 8 + 5\ =\ 6 + 7 \ =\ 7 + 6\\
 &\ =\ 4+4+5\ =\ 4 + 5 + 4 \ =\ 5 + 4 + 4.\end{align*}
Now we define, for $n, p\in \mathbb{N}$,  
\begin{align*}\mathcal{S}^p_n &\ :=\ \{S\subset \{1, \ldots, n\}\,:\, \min S > |S|^p\mbox{ and } n\in S\},\mbox{ and }\\
\mathcal{A}^p_n &\ :=\ \{S\subset \{1, \ldots, n\}\,:\, \min S \geqslant |S|^p\mbox{ and } n\in S\}.
\end{align*}
\begin{theorem}\label{m3} 
For all $n, p\in \mathbb{N}$, it holds that 
\begin{equation}\label{e4}|\mathcal{S}^p_n|\ =\ K_{n, p-1}.\end{equation}
Furthermore, 
$|\mathcal{S}^p_1| = 0$ and $|\mathcal{S}^p_{n+1}| = |\mathcal{A}^p_{n}|$. Hence,
\begin{equation}\label{e3}
    |\mathcal{A}^p_n| \ =\ |\mathcal{S}^p_{n+1}|\ =\ K_{n+1, p-1}.
\end{equation}
\end{theorem}
\begin{remark}\normalfont
When $p=1$, \eqref{e3} and Corollary \ref{c1} give $|\mathcal{A}_n^1| = F_n$. This is Bird's result.  
When $p = 1$, \eqref{e4} and Corollary \ref{c1} give
$|\mathcal{S}^1_n|= F_{n-1}$ for all $n\geqslant 1$. This is precisely the equality ``$A_n = F_{n-1}$" in \cite[Theorem 1]{C1}.
\end{remark}

In Section \ref{linear}, we first go over some preliminary results such as the generalized Pascal triangle and the star-and-bar problem then prove Theorems \ref{m1}, \ref{m2}, and \ref{m4}. In Section \ref{nonlinear}, we study nonlinear Schreier conditions and prove Theorem \ref{m3}.

\section{Schreier conditions on multisets}\label{linear}

Before proving the main results, we briefly mention some preliminary results that will be used in due course. First, for $n\geqslant k\geqslant 0$, the binomial
$\binom{n}{k}$ is located in the $n$\textsuperscript{th} row and $k$\textsuperscript{th} column of the Pascal triangle and represents 
the number of ways of choosing $k$ objects out of $n$ distinguishable objects. In other words, we have $n$ labelled boxes of capacity $1$, and $\binom{n}{k}$ counts the number of ways we can assign $k$ identical objects into these $n$ boxes. This combinatorial interpretation provides the following generalization of the classical binomials and the Pascal triangle. 

Let $s\geqslant 1$ be an integer and let $\binom{n}{k}_s$
denote the number of different ways of distributing $k$ identical objects
among $n$ labelled boxes, each of which may contain at most $s$ objects. Here we require $0\leqslant k\leqslant sn$. Call the triangle whose $n$\textsuperscript{th} row and $k$\textsuperscript{th} column is $\binom{n}{k}_s$ the \textit{$s$-Pascal triangle}. When $s = 1$, we have the Pascal triangle, where the elements in the $n$\textsuperscript{th} row can be constructed by the elements in the $(n-1)$\textsuperscript{th} row through the formula
$$\binom{n}{k} \ =\ \binom{n-1}{k} + \binom{n-1}{k-1},\mbox{ for all } n\geqslant k\geqslant 1.$$
We can do the same for the general $s$-Pascal triangle, using the well-known formula (see \cite[(1.13)]{Bon})
\begin{equation}\label{e2}\binom{n}{k}_{s}\ =\ \sum_{j=0}^{s}\binom{n-1}{k-j}_s.\end{equation}

The following formula in \cite[Theorem 3.1]{BBK} connects $s$-Pascal triangle with the $(s+1)$-step Fibonacci sequence $(F_{n}^{(s+1)})$
\begin{equation}\label{summ}
\sum_{k=0}^{\left\lfloor sn/(s+1)\right\rfloor}\binom{n-k}{k}_s=F_{n+1}^{(s+1)}, \mbox{ for all } n\geqslant 0\mbox{ and } s\geqslant 1.
\end{equation}
For instance, if $s=2$, then \eqref{e2} helps us build the $2$-Pascal triangle 
\begin{center}
	\begin{tabular}{ccccccccc}
		1 &  &  &  &  & & \\
		1 & 1 & 1 &  &  & & \\
		1 & 2 & 3 & 2 & 1 & & \\
		1 & 3 & 6 & 7 & 6 & 3 & 1 &\\
		1 & 4 & 10 & 16 & 19 & 16 & 10 & 4  & 1  \\
          &   &     &   &   \vdots & & & &
	\end{tabular}
\end{center}
and the ascending diagonal sums provide the Tribonacci sequence $(T_n)=(F_n^{(3)})$, the first few terms are
\begin{align*}
	1&\ =\ T_1, \\
	1&\ =\ T_2, \\
	1+1=2&\ =\ T_3,\\
	1+2+1=4&\ =\ T_4, \\
	1+3+3=7&\ =\ T_5, \\
\end{align*}
etc.

Lastly, we recall the star-and-bar lemma. For its proof, see \cite[Lemma 2.1]{KKMW}.
\begin{lemma} \label{l1}
The number of solutions to $x_1 + x_2 + \cdots + x_p = n$ with $x_i\geqslant c_i$ (for some nonnegative number $c_i$) is $\binom{n-\sum_{i=1}^p c_i + (p-1)}{p-1}$.
\end{lemma}

We are now ready to prove our main results.
\begin{proof}[Proof of Theorem \ref{m1}]
	Trivially, $\{n\}\in{\mathcal A}_n^{(s-1)}$. Any extension of the set $\{n\}$ by $k=1,\dots,K_{s-1}$ elements ($K_{s-1}$ to be defined later) is a choice of $k$ elements from the multiset $$\{\underbrace{k+1,\dots,k+1}_{s-1},\dots,\underbrace{n-1,\dots,n-1}_{s-1}\}.$$ 
	Equivalently, we want to put $k$ elements into the boxes labelled by $k+1,k+2,\dots,n-1$, where the capacity of each box is $s-1$. There are $\binom{n-1-k}{k}_{s-1}$ possible choices. 
	In order to determine the precise value of $K_{s-1}$, we consider the inequality 
	$$k\ \leqslant\ (n-k-1)(s-1).$$
	Consequently, $k\leqslant (n-1)(s-1)/s$ and hence, $K_{s-1}=\lfloor(n-1)(s-1)/s\rfloor$. Therefore, we obtain by \eqref{summ} that 
$$
|{\mathcal A}_n^{(s-1)}|\ =\ \sum_{k=0}^{K_{s-1}}\binom{n-1-k}{k}_{s-1}\ =\ F^{(s)}_n.
$$
This completes our proof. 
\end{proof}

\begin{proof}[Proof of Theorem \ref{m2}]
Fix $u\geqslant 2$. For $1\leqslant n\leqslant u$, 
$$\mathcal{B}^{\mathbf s, u}_{n}\ =\ \{\emptyset\},$$
because if $A\in \mathcal{B}^{\mathbf s, u}_{n}$ and $A\neq \emptyset$, then $\min A\geqslant u+1$; however, $$A\ \subset\ \{\underbrace{1, \ldots, 1}_{s_1}, \underbrace{2, \ldots, 2}_{s_2}, \ldots, \underbrace{u, \ldots, u}_{s_u}\}.$$
Hence, $|\mathcal{B}^{\mathbf s, u}_n| = 1$. Choose $n\geqslant u+1$ for $1\leqslant n\leqslant u$. Let us show that 
$$|\mathcal{B}^{\mathbf s, u}_{n}| \ =\ |\mathcal{B}^{\mathbf s, u}_{n-1}| +  |\mathcal{B}^{\mathbf s, u}_{n-u}|.$$
Consider the number of $k$-element sets in $\mathcal{B}^{\mathbf s, u}_{n}$. These sets are subsets of
$$\{\underbrace{uk+1, \ldots, uk+1}_{s_{uk+1}}, \underbrace{uk+2, \ldots, uk+2}_{s_{uk+2}}, \ldots, \underbrace{n, \ldots, n}_{s_n}\}$$
and thus, there are $\binom{k+n-uk-1}{n-uk-1}$ of them according to Lemma \ref{l1}. Call this collection $\mathcal{B}^{\mathbf s, u}_{n,k}$ to have
\begin{equation}\label{e10}\left|\mathcal{B}^{\mathbf s, u}_{n,k}\right|\ =\ \binom{k+n-uk-1}{n-uk-1}\ =\ \binom{k+n-uk-1}{k}.\end{equation}
Here $k\leqslant (n-1)/u$. 
It follows that
\begin{equation}\label{e30}|\mathcal{B}^{\mathbf s, u}_{n}|\ =\ \sum_{k=0}^{\lfloor(n-1)/u\rfloor}\binom{k+n-uk-1}{k},\end{equation}
and in particular, the largest set in $\mathcal{B}_n^{\mathbf s, u}$ has size $\lfloor(n-1)/u\rfloor$. We proceed by case analysis. 

Case 1: $n\not\equiv 1\mod u$. It is easy to check that 
$\lfloor (n-1)/u\rfloor\ =\ \lfloor (n-2)/u\rfloor =: m$. Hence, for any $0\leqslant k\leqslant m$, if $\mathcal{B}^{\mathbf s, u}_n$ has sets of size $k$, then $\mathcal{B}^{\mathbf s, u}_{n-1}$ also has sets of size $k$. By \eqref{e10}, we know that \begin{equation}\label{e11}\mathcal{B}^{\mathbf s, u}_{n-1, k}\ =\ \binom{k+(n-1)-uk-1}{k}\ =\ \binom{n-uk+k-2}{k}, \quad 0\leqslant k\leqslant m.
\end{equation}
Since $\lfloor (n-u-1)/u\rfloor = \lfloor (n-1)/u\rfloor - 1$, $\mathcal{B}^{\mathbf s, u}_{n-u}$ has sets of sizes ranging from $0$ to $m-1$. We have
\begin{align}\label{e12}\mathcal{B}^{\mathbf s, u}_{n-u, k-1}&\ =\ \binom{(k-1)+(n-u)-u(k-1)-1}{k-1}\nonumber\\
&\ =\ \binom{n-uk+k-2}{k-1},\quad 1\leqslant k\leqslant m.
\end{align}
From \eqref{e10}, \eqref{e11}, and \eqref{e12}, we obtain
$$|\mathcal{B}^{\mathbf s, u}_{n, k}|\ =\ |\mathcal{B}^{\mathbf s, u}_{n-1, k}|+|\mathcal{B}^{\mathbf s, u}_{n-u, k-1}|,\quad 1\leqslant k\leqslant m.$$
and thus,
\begin{align*}|\mathcal{B}^{\mathbf s, u}_{n}| - |\mathcal{B}^{\mathbf s, u}_{n, 0}|\ =\ \sum_{k=1}^{m}|\mathcal{B}^{\mathbf s, u}_{n, k}|&\ =\ \sum_{k=1}^{m}|\mathcal{B}^{\mathbf s, u}_{n-1, k}| + \sum_{k=1}^{m}|\mathcal{B}^{\mathbf s, u}_{n-u, k-1}|\\
&\ =\  |\mathcal{B}^{\mathbf s, u}_{n-1}| - |\mathcal{B}^{\mathbf s, u}_{n-1, 0}| + \sum_{k=0}^{m-1}|\mathcal{B}^{\mathbf s, u}_{n-u, k}|.
\end{align*}
Hence, 
$$|\mathcal{B}^{\mathbf s, u}_{n}|\ =\ |\mathcal{B}^{\mathbf s, u}_{n-1}| + |\mathcal{B}^{\mathbf s, u}_{n-u}|.$$

Case 2: $n\equiv 1\mod u$. Let $m  = (n-1)/u$. Then $\lfloor (n-2)/u\rfloor = m-1$ and so, $\mathcal{B}^{\mathbf s, u}_{n}$ does not contain any set of size $m$. By \eqref{e10}, 
\begin{equation}\label{e15}\mathcal{B}^{\mathbf s, u}_{n-1, k}\ =\ \binom{k+(n-1)-uk-1}{k}\ =\ \binom{n-uk+k-2}{k}, \quad 0\leqslant k\leqslant m-1.
\end{equation}
Similarly, sets in $\mathcal{B}$ have sizes ranging from $0$ to $\lfloor (n-u-1)/u\rfloor = m-1$. Particularly,
\begin{equation}\label{e16}\mathcal{B}^{\mathbf s, u}_{n-1, k-1}\ =\ \binom{n-uk+k-2}{k-1}, \quad 1\leqslant k\leqslant m.\end{equation}
From \eqref{e10}, \eqref{e15}, and \eqref{e16}, we obtain 
$$|\mathcal{B}_{n,k}^{\mathbf s, u}|\ =\ |\mathcal{B}_{n-1,k}^{\mathbf s, u}| + |\mathcal{B}^{\mathbf s, u}_{n-u,k-1}|, \quad 1\leqslant k\leqslant m-1.$$
Hence,
$$\sum_{k=1}^{m-1}|\mathcal{B}^{\mathbf s, u}_{n,k}| \ =\ \sum_{k=1}^{m-1}|\mathcal{B}^{\mathbf s, u}_{n-1,k}| + \sum_{k=0}^{m-2}|\mathcal{B}^{\mathbf s, u}_{n-u,k}|,$$
which gives
$$|\mathcal{B}_n^{\mathbf s, u}| - |\mathcal{B}_{n,0}^{\mathbf s, u}| - |\mathcal{B}^{\mathbf s, u}_{n,m}|\ =\ \left(|\mathcal{B}_{n-1}^{\mathbf s, u}| - |\mathcal{B}^{\mathbf s, u}_{n-1,0}|\right) + (|\mathcal{B}_{n-u}^{\mathbf s, u}| - |\mathcal{B}^{\mathbf s, u}_{n-u, m-1}|).$$
Since $|\mathcal{B}^{\mathbf s, u}_{n,0}| = |\mathcal{B}^{\mathbf s, u}_{n-1,0}|$, it remains to verify that
$|\mathcal{B}^{s,u}_{n,m}| = |\mathcal{B}^{\mathbf s, u}_{n-u, m-1}|$, which holds due to \eqref{e10} and $n-1 = um$. This completes our proof. 
\end{proof}

\begin{proof}[Proof of Theorem \ref{m4}]
Fix $n\in\mathbb{N}$ and $s\geqslant 2$. Let $C\in \mathcal{C}^{(s-1)}_n$ and $|C|=k\geqslant 1$ for some $1\leqslant k\leqslant n$. Then we must choose $k-1$ elements from the set
$$
\{k_1,k_2,\dots,k_{s-1},\dots,(n-1)_1,(n-2)_2,\dots,(n-1)_{s-1}\}\ \subset\ H_n.
$$	
The number of ways to do so is obviously $\binom{(n-k)(s-1)}{k-1}$. Here $k\leqslant \lfloor (n(s-1)+1)/s\rfloor$. Hence,
\begin{align*}
|\mathcal{C}^{(s-1)}_n|&\ =\ \sum_{k=1}^{\lfloor (n(s-1)+1)/s\rfloor}\binom{n(s-1)-ks+k}{k-1}\\
&\ =\ \sum_{k=0}^{\lfloor (n-1)(s-1)/s\rfloor}\binom{(n-1)(s-1)-ks+k}{k}\ =\ |\mathcal{B}^{\mathbf s, s}_{(n-1)(s-1)+1}|,
\end{align*}
where the last equality is due to \eqref{e30}. By Theorem \ref{m2}, we conclude that
$$|\mathcal{C}^{(s-1)}_n|\ =\ K^{(s)}_{(n-1)(s-1)+1}.$$
This completes our proof.
\end{proof}

\section{Nonlinear Schreier conditions}\label{nonlinear}
To find a formula for $K_{n,p}$, we again use the star-and-bar problem.
\begin{proposition}\label{p1}
Fix $n\in\mathbb{N}$ and $p\in \mathbb{Z}_{\geqslant 0}$. Let $K_{n,p,k}$ be the number of decompositions of $n$ into $k$ parts such that the smallest part is strictly greater than $k^p$. It holds that
$$K_{n, p,k} \ = \ \binom{n-k^{p+1}-1}{k-1}.$$
Hence, 
\begin{equation}\label{e1}K_{n,p}\ =\ \sum_{k=1}^n K_{n,p,k}\ =\ \sum_{k=1}^{n}\binom{n-k^{p+1}-1}{k-1}\end{equation}
with the convention that $\binom{u}{v} = 0$ if $u < v$.
\end{proposition}
\begin{proof}
Clearly, $K_{n,p,k}$ is equal to the number of solutions to $x_1 + \cdots + x_k = n$ such that $x_i\geqslant k^p+1$ for all $1\leqslant i\leqslant k$. By Lemma \ref{l1}, 
$$K_{n,p,k}\ =\ \binom{n+(k-1)-k(k^p+1)}{k-1}\ =\ \binom{n-k^{p+1}-1}{k-1}.$$

The second statement follows directly from the definition of $K_{n,p, k}$.
\end{proof}
\begin{corollary}\label{c1}
We have 
$$K_{n, 0} \ =\ F_{n-1}, \quad n\in \mathbb{N}.$$
\end{corollary}
\begin{proof}
    By Proposition \ref{p1}, 
    $$K_{n,0}\ =\ \sum_{k=1}^{n}\binom{n-k-1}{k-1}\ =\ \sum_{k=0}^{n-1}\binom{n-2-k}{k}\ =\ F_{n-1}, \quad n\in \mathbb{N}.$$
\end{proof}

\begin{proof}[Proof of Theorem \ref{m3}]
We calculate $|\mathcal{S}^p_n|$ by considering sets $S\in \mathcal{S}^p_n$ of certain size $k\geqslant 1$. Since $S\in \mathcal{S}^p_n$, we have $\min S\geqslant k^p+1$. Therefore, $S\backslash \{n\}\subset \{k^p+1, k^p+2, \ldots, n-1\}$. Since $|S\backslash \{n\}| = k-1$, it follows that the number of choices for $S\backslash \{n\}$ is $\binom{n-k^p-1}{k-1}$. Hence, 
\begin{equation}\label{e20}|\mathcal{S}^p_n|\ =\ \sum_{k=1}^{n}\binom{n-k^p-1}{k-1}\ =\ K_{n, p-1},\end{equation}
due to Proposition \ref{p1}. 

For the second statement, we clearly have $|\mathcal{S}^p_1| = 0$. Using the same reasoning as above and \eqref{e20}, we have 
$$|\mathcal{A}^p_n|\ =\ \sum_{k=1}^n \binom{n-k^p}{k-1}\ =\ \sum_{k=1}^{n+1}\binom{(n+1)-k^p-1}{k-1}\ =\ |\mathcal{S}^p_{n+1}|.$$
\end{proof}

Finally, we modify $\mathcal{A}_n^p$ so that the maximum of sets is not fixed. Fix $p\in \mathbb{N}$ and define 
$$\mathcal{B}^p_n \ :=\ \{S\subset \{1,\ldots, n\}\,:\, \mbox{ either }S = \emptyset \mbox{ or } (\max S - 1\in S\mbox{ and }\min S \geqslant |S|^p)\}.$$
\begin{theorem}
    It holds that $$|\mathcal{A}^p_n| \ =\ |\mathcal{B}^p_n|, \quad n, p\in \mathbb{N}.$$
\end{theorem}
\begin{proof} Fix $n, p\in\mathbb{N}$.
    For $k\geqslant 1$, let 
    $$\mathcal{A}^p_{n, k} \ =\ \{S\in \mathcal{A}^p_n\,:\, |S| = k\}\mbox{ and }\mathcal{B}^p_{n, k} \ =\ \{S\in \mathcal{B}^p_n\,:\, |S| = k\}.$$
    Observe that $|\mathcal{A}^p_{n,1}| = 1$, while $|\mathcal{B}^p_{n,1}| = 0$. However, this discrepancy is remedied by the appearance of the empty set in $\mathcal{B}^p_n$. Hence, it suffices to show that
    $$|\mathcal{A}^p_{n,k}| \ =\ |\mathcal{B}^p_{n,k}|, \quad k\geqslant 2.$$
    As in the proof of Theorem \ref{m3}, we have
    \begin{equation}\label{e5}|\mathcal{A}^p_{n,k}|\ =\ \binom{n-k^p}{k-1}.\end{equation}
    On the other hand, if $S\in \mathcal{B}^p_{n,k}$, then $\min S\geqslant k^p$ and so, $S\subset \{k^p, k^p+1, \ldots, n\}$. Since $\max S-1\in S$, the maximum element $m$ of $S$ belongs to $\{k^p+1, \ldots, n\}$. After the maximum is chosen, the second largest element (that is, $m-1$) is fixed. There are $\binom{m-k^p-1}{k-2}$ choices to choose the remaining $k-2$ elements from $\{k^p, \ldots, m-2\}$. Therefore, 
    \begin{equation}\label{e6}|\mathcal{B}^p_{n,k}|\ =\ \sum_{m=k^p+1}^{n}\binom{m-k^p-1}{k-2}\ =\ \sum_{m=k^p+k-1}^n \binom{m-k^p-1}{k-2}\ =\ \binom{n-k^p}{k-1},\end{equation}
    where the last equality is the well-known hockey-stick identity (see \cite[Theorem 1.2.3 item (5)]{W}.) From \eqref{e5} and \eqref{e6}, we have the desired conclusion. 
\end{proof}

We end by mentioning two directions for further research:

First, the linear Schreier-type condition has been thoroughly investigated in \cite{BCF, CMX}. However, less is known about nonlinear Schreier conditions. In this paper, we count sets $F$ that satisfy the nonlinear condition $\min F \geqslant |F|^s$, where $s\in \mathbb{N}_{\geqslant 2}$ and connect their counts to decompositions of integers. Further research can investigate other nonlinear conditions. 

Second, \cite{C2} showed a way to use the Schreier condition to obtain partial sums of an arbitrary order of the Fibonacci sequence (also called the hyperfibonacci sequences.)
Can we obtain partial sums of an arbitrary order of the $s$-step Fibonacci sequence from Schreier-type conditions on multisets? Here, by an arbitrary order $k\geqslant 1$, we mean the sequence obtained by applying the partial sum operator to the $s$-step Fibonacci sequence $k$ times. 

\noindent {\bf Acknowledgement.} This work was completed as part of the 2022 Polymath Jr program. We thank our colleagues there for helpful conversations. For L.~Szalay, the research was supported by National Research, Development and Innovation Office Grant 2019-2.1.11-T\'ET-2020-00165,
by Hungarian National Foundation for Scientific Research Grant No.~128088, and No.~130909, and by the Slovak Scientific Grant Agency VEGA 1/0776/21.

\ \\
\end{document}